\RequirePackage{fix-cm}

\documentclass[smallcondensed]{svjour3}

\newif\ifTikz
\Tikztrue

\smartqed

\usepackage[linesnumbered,vlined,boxed]{algorithm2e}
\usepackage{graphicx}
\usepackage{amsmath,amsfonts}
\usepackage{tensor}
\usepackage{subfig}
\usepackage{pgf,tikz}
\usepackage{pgfplots}
\usepackage{mathtools}
\usepackage[colorlinks=true,citecolor=blue,urlcolor=blue,linkcolor=black]{hyperref}
\usepackage{todonotes}
\usepackage[displaymath, mathlines,running]{lineno}
\usepackage{mathdots}
\usepackage{afterpage}

\DeclareMathOperator\erfc{erfc}

\pgfplotsset{compat=newest}
\usetikzlibrary{spy}
\usetikzlibrary{backgrounds}
\usetikzlibrary{decorations.markings}
\usetikzlibrary{shapes,arrows}

\usepackage{rotating}
\usepackage{pdflscape}

\usepackage{array}
\newcolumntype{x}[1]{>{\centering\arraybackslash\hspace{0pt}}p{#1}}

\begin{document}

\title{Pipeline Implementations of Neumann--Neumann and
  Dirichlet--Neumann Waveform Relaxation Methods}

\titlerunning{Pipeline NNWR and DNWR Methods}
\author{
  Benjamin W.~Ong \and
  Bankim C.~Mandal
}
\authorrunning{B.~W.~Ong, B.~C.~Mandal}

\institute{
  Benjamin W.~Ong \at
  Mathematical Sciences, Michigan Technological University, Houghton, MI, 49931\\
  {\tt ongbw@mtu.edu}
  \and
  Bankim C.~Mandal \at Dept.~of
  Mathematics, Michigan State University, East Lansing, MI, 48824 \\
  {\tt bmandal@math.msu.edu}
}

\date{Received: \today}
         
\maketitle

\begin{abstract}
  This paper is concerned with the reformulation of Neumann-Neumann
  Waveform Relaxation (NNWR) methods and Dirichlet-Neumann Waveform
  Relaxation (DNWR) methods, a family of parallel space-time
  approaches to solving time-dependent PDEs.  By changing the order of
  the operations, pipeline-parallel computation of the waveform
  iterates are possible without changing the final solution.  The
  parallel efficiency and the increased communication cost of the
  pipeline implementation is presented, along with weak scaling
  studies to show the effectiveness of the pipeline NNWR and DNWR
  algorithms.
\end{abstract}

\maketitle


\section{Introduction}

Dirichlet--Neumann waveform relaxation (DNWR) and Neumann--Neumann
waveform relaxation (NNWR) methods \cite{GKM1,GKM2,Mandal,Kwok,GKM3}
have been formulated and analyzed recently for parabolic and
hyperbolic partial differential equations (PDEs).  These iterative
methods are based on non-overlapping domain decomposition in space,
where the iterations require subdomain solves with Dirichlet or Neumann
boundary conditions.  The superlinear convergence behavior of both the
DNWR and NNWR methods have previously been shown for the heat equation
\cite{GKM1,GKM2}; finite step convergence to the exact solution for
the wave equation has also been shown \cite{GKM3,Mandal2}.

Both the DNWR and NNWR methods belong to the family of Waveform
Relaxation (WR) methods.  In the classical WR formulation, one solves
the space-time subproblem over the entire time horizon at each
iteration, before communicating interface data across subdomains.
These WR methods have their origin in the work of Picard--Lindel\"of
\cite{Picard,Lind} in the 19th century; WR methods were introduced as
a parallel approach for solving systems of ODEs \cite{LelRue}. Since
then, the WR framework has been combined with numerical approaches for
solving elliptic problems for tackling time-dependent parabolic PDEs
\cite{MR1638096,GilKel}. New transmission conditions have also been
developed, known as optimized Schwarz WR (OSWR) methods, in order to
achieve convergence with non-overlapping subdomains, or in general,
faster convergence for overlapping subdomains; OSWR approaches for
parabolic problems \cite{MR2300292} and OSWR approaches for hyperbolic
problems \cite{MR2035001} have both been explored.  This paper discusses a
pipeline implementation of the DNWR and NNWR algorithms
\cite{BankThes}, which are extensions of Dirichlet-Neumann
\cite{MR865945,MR829613,MR998911,MarQuar02} and Neumann-Neumann
\cite{MR992000,MR1106455,MR1095198,MR2104179} algorithms for solving space--time PDEs.  By
carefully re-arranging the order of operations, different waveform
iterates can be simultaneously computed in a pipeline-parallel
fashion, {\em without} changing the final solution.

The pipeline implementation of the DNWR and NNWR algorithm subdivides
the entire time domain into smaller time blocks.  Each space-time
subproblem is solved on this smaller time block, and updated
transmission conditions are transmitted before advancing to the next
time block.  This enables many concurrent subdomain space-time blocks,
resulting in efficient parallelization on modern
supercomputers. Specifically, given an appropriate number of
processors, many waveform iterates can be computed in the same
wall-clock time as a single processor computing one waveform iterate.
The pipeline implementation for WR relaxation was mentioned in
\cite{MR1340665,MR1146977}, and
recently reintroduced and benchmarked for the classical Schwarz WR
method \cite{pswr-dd22}.

The theoretical presentation of the algorithms described in this paper
are for one-dimensional time-dependent PDEs of the form
\begin{subequations}
  \label{eqn:pde}
  \begin{align}
    \partial_t u - \mathcal{L}u  =  f(x,t), 
    \quad (x,t)\in\Omega, \\
    u(x,0)  =  u_{0}(x), 
    \quad x\in [0,L],\\
    u(0,t)  =  g_l(t),
    \quad
    u(L,t) = g_r(t),
    \quad t \in [0,T],
  \end{align}
\end{subequations}
where $\mathcal{L}$ is a spatial operator, the space-time domain,
$\Omega:[0,L]\times[0,T]$, is a bounded domain, $\partial\Omega$ is a
smooth boundary, $u_0(x)$ is the initial condition, and $g_l(t)$ and
$g_r(t)$ are Dirichlet boundary conditions.  The pipeline
implementations can also be applied naturally to time-dependent PDEs
of the form $\partial_{tt} - \mathcal{L}(u)=f(x,t)$, as well as WR
methods formulated for solving PDEs in higher spatial dimensions.

This paper is broken into two main sections. In
Section~\ref{sec:NNWR}, we review the NNWR method for solving equation
\eqref{eqn:pde} before presenting the pipeline approach and
numerical studies.  The DNWR algorithm is outlined in
Section~\ref{sec:DNWR}, along with various pipeline approaches and
numerical results.

\section{Neumann--Neumann Waveform Relaxation (NNWR)}
\label{sec:NNWR}

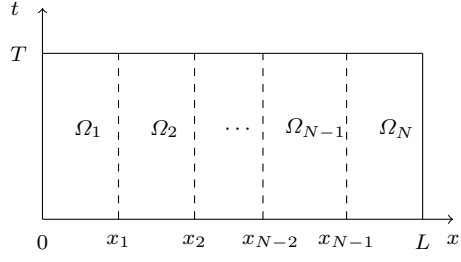
\begin{figure}
  \centering
  \begin{tikzpicture}
 
    \draw[->] (0,2.2) -- (5.4,2.2);
    \draw[->] (0,2.2) -- (0,5); 
    \draw[-] (0,4.4) -- (5,4.4);
    \draw[-] (5,2.2) -- (5,4.4);
    
    \draw[dashed] (1,2.2) -- (1,4.4);
    \draw[dashed] (2,2.2) -- (2,4.4);
    \draw[dashed] (2.9,2.2) -- (2.9,4.4);
    \draw[dashed] (4,2.2) -- (4,4.4);

    \draw (0.9,3.4) node[left] {$\Omega_1$};
    \draw (1.9,3.4) node[left] {$\Omega_2$};
    \draw (4.1,3.4) node[left] {$\Omega_{N-1}$};
    \draw (2.9,3.4) node[left] {$\cdots$};
    \draw (5,3.4) node[left] {$\Omega_N$};
    \draw (0,2.1) node[below] {$0$};
    \draw (5,2.1) node[below] {$L$};
    \draw (-0.2,5) node[left] {$t$};
    \draw (-0.1,4.4) node[left] {$T$};
    \draw (5.4,2.1) node[below] {$x$};
    
    \draw (1,2.1) node[below] {$x_1$};
    \draw (2,2.1) node[below] {$x_2$};
    \draw (3,2.1) node[below] {$x_{N-2}$};
    \draw (4,2.1) node[below] {$x_{N-1}$};    
  \end{tikzpicture}
  
  \caption{Decomposition of a space-time domain, $[0,L]\times[0,T]$ domain} 
  \label{fig:subdomains}
\end{figure}

The NNWR algorithm for the model problem \eqref{eqn:pde} with
multi--subdomains setting was previously proposed and analyzed in
\cite{GKM1}. In this method, the space-time domain $\Omega:[0,L]\times[0,T]$
is first partitioned into non-overlapping space-time subdomains, $\{\Omega_{i}$,
$1\leq i\leq N\}$. Figure~\ref{fig:subdomains} illustrates a simple 1D
spatial decomposition, although more complicated 2D or 3D
decompositions (with cross points) are also possible.
Let $\partial\Omega_i$ denote the domain boundary of $\Omega_i$, and
let $w_{i}^{[0]} = w(x_i,t), i =1,\ldots,N-1$, be an initial
time-dependent guess on the subdomain boundaries.  The NNWR method
performs a two-step iteration for each waveform iterate,
$k=1,2,\ldots$, until convergence is reached.  This two-step iteration
consists of first solving a ``Dirichlet subproblem'' on each
space-time subdomain,
\begin{subequations}
  \label{eqn:nnwr_dirichlet}
  \begin{align}
    \partial_{t}u_{i}^{[k]} - \mathcal{L}u_{i}^{[k]}  =  f,
    &\quad (x,t) \in \Omega_i,\\
    u_{i}^{[k]}(x,0)  =  u_{0}(x),
    &\quad x \in (x_{i-1},x_i),\\
    u_{i}^{[k]}(x_{i-1},t) =
    &\begin{dcases}
       g_l(t) & \text{if } i = 1 \\
      w_{i-1}^{[k-1]}(t) & \text{otherwise}
    \end{dcases}, \\
    u_{i}^{[k]}(x_{i},t) =
    &\begin{dcases}
      g_r(t) & \text{if } i = N \\
      w_{i}^{[k-1]}(t) & \text{otherwise}
    \end{dcases}, 
  \end{align}
\end{subequations}
followed by solving an auxiliary ``Neumann subproblem''
\begin{subequations}
  \label{eqn:nnwr_neumann}
  \begin{align}
    &\partial_{t}\psi_{i}^{[k]}-\mathcal{L}\psi_{i}^{[k]}  =  0,
    \quad (x,t) \in \Omega_i\\
    &\psi_{i}^{[k]}(x,0)  =  0,
    \quad x \in (x_{i-1},x_i),\\
    &\begin{cases}
    \psi_1^{[k]}(0) = 0,  &\text{if } i = 1,\\
      -\partial_{x}\psi_{i}^{[k]}(x_{i-1},t) =
      (\partial_{x}u_{i-1}^{[k]}-\partial_{x}u_{i}^{[k]})(x_{i-1},t),  &\text{if } i>1,
     \end{cases}\\
    &\begin{cases}
       \partial_{x}\psi_{i}^{[k]}(x_{i},t) =
       (\partial_{x}u_{i}^{[k]}-\partial_{x}u_{i+1}^{[k]})(x_{i},t),  &\text{for } i < N, \\
       \psi_N^{[k]}(L) = 0,  &\text{if } i = N.
     \end{cases}
  \end{align}
\end{subequations}
Then, the Dirichlet traces at the subdomain interfaces  are updated,
\begin{align}
  \label{eqn:nnwr_dirichlet_update}
  w_{i}^{[k]}(t)=w_{i}^{[k-1]}(t)-\theta
  \left(\psi_{i}^{[k]}(x_{i},t)+\psi_{i+1}^{[k]}(x_{i},t)\right).
\end{align}
\begin{remark}
  The auxiliary equations are only solved if the waveform iterates
  have not converged.  The convergence of the NNWR algorithm for the
  heat equation was previously analyzed \cite{GKM1}.  For
  $\theta=1/4$, the NNWR algorithm converges superlinearly with the
  estimate
  \begin{align*}
    \displaystyle
    \max_{1\leq i\leq N-1}\| w_{i}^{[k]}\|_{L^{\infty}(0,T)}
     \leq\left(\frac{\sqrt{6}}{1-e^{-\frac{(2k+1)\tilde{h}^{2}}{T}}}\right)^{2k}
    e^{-k^{2}\tilde{h}^{2}/T}\max_{1\leq i\leq N-1}
    \| w_{i}^{[0]}\|_{L^{\infty}(0,T)},
  \end{align*}
  where $\tilde{h}=\min\{x_i-x_{i-1} : i=1,\ldots,N\}$ is the minimum
  subdomain width.

  For the wave equation, the NNWR algorithm converges after the second
  waveform iterate, provided the window of integration is sufficiently
  small.
\end{remark}
  
\subsection{Classical NNWR implementation}

A pseudo-code for a classical implementation of the NNWR in
$\mathbb{R}^1$ is given in Algorithm~\ref{alg:NNWR}.  If the domain is
broken into $N$ non-overlapping subdomains, the classical
implementation uses a total of $N$ processing cores to approximate the
solution. $2(N-1)(2K-1)$ or $4(N-1)K$ messages are needed (depending
on when the stopping criterion is satisfied), with each message
containing $N_t$ words, where $K$ is the number of waveform iterates
computed and $N_t$ is the number of time steps.  In this simplified
pseudo-code, it is assumed that identical time discretizations are
taken in each subdomain.  A generalization where each subdomain might
take different discretizations is possible, provided an interpolation
of the Neumann or Dirichlet traces are computed before lines 10 and 16
of Algorithm~\ref{alg:NNWR}.
\begin{algorithm}
  \ForPar{$i=1$ \KwTo $N$}{
    \If {$i < N$} {
      Guess $w_i^{[0]}(t_\ell), \quad \ell =1,\ldots,N_t$\;
    }
    \If {$i > 1 $} {
      Guess $w_{i-1}^{[0]}(t_\ell),\quad \ell = 1,\ldots,N_t$\;
    }
    Set $k=1$\;
    \While{not converged}{
      \For{$\ell=1$ \KwTo $N_t$}{
        Solve equation~\eqref{eqn:nnwr_dirichlet} for $u_i^{[k]}(x,t_\ell)$\;
        Compute jump in Neumann data along interfaces\;
      }
      Transmit/Receive Neumann data\;
      Check for convergence.  If converged, {\tt break}\;
      \If{not converged}{
        \For{$\ell = 1$ \KwTo $N_t$}{
          Solve equation~\eqref{eqn:nnwr_neumann} for $\psi_i^{[k]}(x,t_\ell)$\;
          Solve equation~\eqref{eqn:nnwr_dirichlet_update}  for $w_i^{[k]}(t_\ell)$\;
        }
        Transmit/Receive Dirichlet data\;
        Check for convergence\;
      }
      $k \leftarrow k+1$\;
    }
  }
\caption{The classical implementation of the NNWR algorithm is able to
  utilize $N$ computing cores if the domain is broken into $N$
  non-overlapping subdomains.}
\label{alg:NNWR}  
\end{algorithm}

\subsection{Pipeline NNWR implementation}
\label{sec:pipeline_NNWR}

\begin{figure}
  \centering
  \begin{tikzpicture}
  
  \tikzstyle{p1} = [rectangle, draw, fill=blue!20, rounded corners]
  \tikzstyle{p2} = [rectangle, draw, fill=orange!20, rounded corners]
  \tikzstyle{p3} = [rectangle, draw, fill=cyan!20, rounded corners]
  \tikzstyle{p4} = [rectangle, draw, fill=red!20, rounded corners]

  \node[p1] (u1w11) {$u^{[1]} \in \Omega_{11}$};
  \node[p1, right = of u1w11] (u1w12) {$u^{[1]} \in \Omega_{12}$};
  \node[p3, below = of u1w12] (p1w11) {$\psi^{[1]} \in \Omega_{11}$};
  \node[p1, right = of u1w12] (u2w11) {$u^{[2]} \in \Omega_{11}$};
  \node[p3, below = of u2w11] (p1w12) {$\psi^{[1]} \in \Omega_{12}$};
  \node[p1, right = of u2w11] (u2w12) {$u^{[2]} \in \Omega_{12}$};
  \node[p3, below = of u2w12] (p2w11) {$\psi^{[2]} \in \Omega_{11}$};
  \node[p3, right = of p2w11] (p2w12) {$\psi^{[2]} \in \Omega_{12}$};

  \node[p2, below = of p1w11] (u1w22) {$u^{[1]} \in \Omega_{22}$};
  \node[p2, left = of u1w22] (u1w21) {$u^{[1]} \in \Omega_{21}$};
  \node[p4, below = of u1w22] (p1w21) {$\psi^{[1]} \in \Omega_{21}$};
  \node[p2, right = of u1w22] (u2w21) {$u^{[2]} \in \Omega_{21}$};
  \node[p4, below = of u2w21] (p1w22) {$\psi^{[1]} \in \Omega_{22}$};
  \node[p2, right = of u2w21] (u2w22) {$u^{[2]} \in \Omega_{22}$};
  \node[p4, below = of u2w22] (p2w21) {$\psi^{[2]} \in \Omega_{21}$};
  \node[p4, right = of p2w21] (p2w22) {$\psi^{[2]} \in \Omega_{22}$};

  \begin{scope}[color=blue!50,decoration={
        markings,
        mark=at position 0.4 with {\arrow{>}}}
      ] 
    \draw[postaction={decorate}] (u1w11.east) -- (p1w11.west);
    \draw[postaction={decorate}] (u1w11.east) -- (p1w21.west);
    \draw[postaction={decorate}] (u1w12.east) -- (p1w12.west);
    \draw[postaction={decorate}] (u1w12.east) -- (p1w22.west);
    \draw[postaction={decorate}] (u2w11.east) -- (p2w11.west);
    \draw[postaction={decorate}] (u2w11.east) -- (p2w21.west);
    \draw[postaction={decorate}] (u2w12.east) -- (p2w12.west);
    \draw[postaction={decorate}] (u2w12.east) -- (p2w22.west);
  \end{scope}
  
  \begin{scope}[color=orange!50,decoration={
        markings,
        mark=at position 0.4 with {\arrow{>}}}
      ] 
    \draw[postaction={decorate}] (u1w21.east) -- (p1w11.west);
    \draw[postaction={decorate}] (u1w21.east) -- (p1w21.west);
    \draw[postaction={decorate}] (u1w22.east) -- (p1w12.west);
    \draw[postaction={decorate}] (u1w22.east) -- (p1w22.west);
    \draw[postaction={decorate}] (u2w21.east) -- (p2w11.west);
    \draw[postaction={decorate}] (u2w21.east) -- (p2w21.west);
    \draw[postaction={decorate}] (u2w22.east) -- (p2w12.west);
    \draw[postaction={decorate}] (u2w22.east) -- (p2w22.west);
  \end{scope}

  \begin{scope}[color=cyan!50,decoration={
        markings,
        mark=at position 0.4 with {\arrow{>}}}
    ] 
    \draw[postaction={decorate}] (p1w11.east) -- (u2w11.west);
    \draw[postaction={decorate}] (p1w11.east) -- (u2w21.west);
    \draw[postaction={decorate}] (p1w12.east) -- (u2w12.west);
    \draw[postaction={decorate}] (p1w12.east) -- (u2w22.west);
  \end{scope}

  \begin{scope}[color=red!50,decoration={
        markings,
        mark=at position 0.4 with {\arrow{>}}}
    ] 
    \draw[postaction={decorate}] (p1w21.east) -- (u2w11.west);
    \draw[postaction={decorate}] (p1w21.east) -- (u2w21.west);
    \draw[postaction={decorate}] (p1w22.east) -- (u2w12.west);
    \draw[postaction={decorate}] (p1w22.east) -- (u2w22.west);
  \end{scope}
  
\end{tikzpicture}
  \caption{A graphical illustration of a pipeline NNWR implementation
    for a two-domain, two time block example.  The different color
    blocks represent different processing cores -- in this case, four
    processors can be used to compute this pipeline NNWR
    implementation.}
  \label{fig:pipeline_NNWR}
\end{figure}
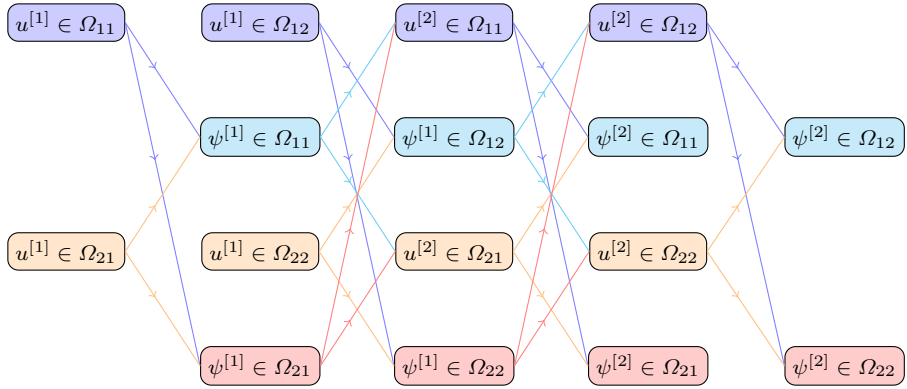

A pipeline implementation of the NNWR algorithm allows for higher
concurrency at the expense of increased communication: multiple
waveform iterates are simultaneously evaluated after an initial
start-up cost.  Both implementations will return the same numerical
solution.  The main idea is to decompose the time window into
non-overlapping blocks, and transmit messages to available processors
after the evaluation of each time block, so that other processors can
simultaneously compute solutions to the auxiliary equations or the
next waveform iterate.  We first illustrate this for a simple two
subdomains, two time blocks example.  Let $0=T_0 < T_1 < T_2=T$ and
let $\Omega_{ij}$ denotes the space-time domain, $[x_{i-1},x_i] \times
[T_{j-1},T_j]$.
The pipeline implementation starts with solving the
Dirichlet subproblem, equations~\eqref{eqn:nnwr_dirichlet}, for
$u^{[1]}$ in $\Omega_{11}$ and $\Omega_{21}$ in parallel with two
processors. As soon as this computation is completed, messages are
sent to available processors. The algorithm then proceeds with solving
the Dirichlet subproblems for $u^{[1]}$ in $\Omega_{12}$ and
$\Omega_{22}$, as well as solving auxiliary
equation~\eqref{eqn:nnwr_neumann} for $\psi^{[1]}$ in $\Omega_{11}$
and $\Omega_{21}$ simultaneously, using a total of four processors.
If the waveform iterates have not converged, messages are sent to the
appropriate processors, and the algorithm proceeds with solving the
Dirichlet subproblems for $u^{[2]}$ in $\Omega_{11}$ and
$\Omega_{21}$, as well as solving the auxiliary subproblem for
$\psi^{[1]}$ in $\Omega_{12}$ and $\Omega_{22}$, again with four
processors.  A graphical illustration of the pipeline framework for
this example, assuming two full waveform iterates are desired, is
shown in Figure~\ref{fig:pipeline_NNWR}.

\begin{figure}
  \centering
  \begin{tikzpicture}
 
    \draw[->] (0,2.2) -- (5.4,2.2);
    \draw[->] (0,2.2) -- (0,5); 
    \draw[-] (0,4.6) -- (5,4.6);
    \draw[-] (5,2.2) -- (5,4.6);
    
    \draw[dashed] (1,2.2) -- (1,4.6);
    \draw[dashed] (2,2.2) -- (2,4.6);
    \draw[dashed] (2.8,2.2) -- (2.8,4.6);
    \draw[dashed] (4,2.2) -- (4,4.6);
    
    \draw[dashed](0,2.8) -- (5,2.8);
    \draw[dashed](0,3.4) -- (5,3.4);
    \draw[dashed](0,4) -- (5,4);
    
    \draw (0.9,2.5) node[left] {$\Omega_{11}$};
    \draw (1.9,2.5) node[left] {$\Omega_{21}$};
    \draw (4.1,2.5) node[left] {$\Omega_{N-1,1}$};
    \draw (2.8,2.5) node[left] {$\cdots$};
    \draw (5,2.5) node[left] {$\Omega_{N1}$};
    
    \draw (0.9,3.1) node[left] {$\Omega_{12}$};
    \draw (1.9,3.1) node[left] {$\Omega_{22}$};
    \draw (4.1,3.1) node[left] {$\Omega_{N-1,2}$};
    \draw (2.8,3.1) node[left] {$\cdots$};
    \draw (5,3.1) node[left] {$\Omega_{N2}$};
    
    \draw (0.7,3.8) node[left] {$\vdots$};
    \draw (1.7,3.8) node[left] {$\vdots$};
    \draw (3.6,3.8) node[left] {$\vdots$};
    \draw (2.6,3.8) node[left] {$\vdots$};
    \draw (4.7,3.8) node[left] {$\vdots$};
    
    \draw (0.9,4.3) node[left] {$\Omega_{1J}$};
    \draw (1.9,4.3) node[left] {$\Omega_{2J}$};
    \draw (4.1,4.3) node[left] {$\Omega_{N-1,J}$};
    \draw (2.8,4.3) node[left] {$\cdots$};
    \draw (5,4.3) node[left] {$\Omega_{NJ}$};
    
    \draw (0,2.1) node[below] {$0$};
    \draw (5,2.1) node[below] {$L$};
    \draw (-0.2,5) node[left] {$t$};
    \draw (-0.1,4.4) node[left] {$T$};
    \draw (5.4,2.1) node[below] {$x$};
    
    \draw (1,2.1) node[below] {$x_1$};
    \draw (2,2.1) node[below] {$x_2$};
    \draw (3,2.1) node[below] {$x_{N-2}$};
    \draw (4,2.1) node[below] {$x_{N-1}$};    
  \end{tikzpicture}

  \caption{Decomposition of a space-time domain, $[0,L]\times[0,T]$,
    into $N\times J$ subdomains for pipeline parallelism.}
  \label{fig:pipeline_domains}
\end{figure}
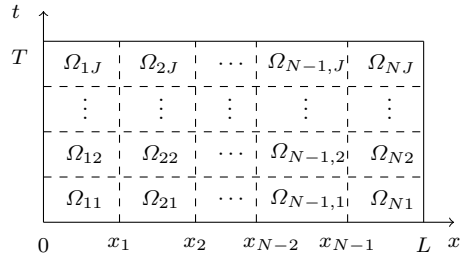

A pseudo-code for a pipeline implementation of the NNWR in
$\mathbb{R}^1$ for the case $J\ge2K$ is given in
Algorithm~\ref{alg:pipeline_NNWR}.  If the domain is broken into
$N\times J$ non-overlapping subdomains as shown in Figure~\ref{fig:pipeline_domains},
this implementation is able to utilize up to $2NK$ processing cores,
where $K$ is the number of waveform iterates computed and $J$ is the
number of time blocks.  The number of messages increases, with up to
$4J(N-1)K$ messages needed in the algorithm.  The size of each message
decreases to $N_t/J$ words.  If $J<2K$, the pipeline implementation is
only able to utilize $NJ$ processors.  A pseudo-code for a pipeline
NNWR implementation for the case $J<2K$ is given in
Algorithm~\ref{alg:pipeline_NNWR_case2} in the Appendix.
\begin{algorithm}
  \tcc{Assumes $J\ge2K$}
  \ForPar{$i=1$ \KwTo $N$}{
    \ForPar{$k=1$ \KwTo $K$}{
      \ForPar{$q=0$ \KwTo $1$}{
        \If{$q==1$}{
          \tcc{Dirichlet Update}
          \For{$j=1$ \KwTo $J$}{
            \eIf{$ k >1$}{
              receive updated Dirichlet data\;
            }{
              \For{$\ell=1$ \KwTo $N_t/J$}{
                $t_\ell = \left(\ell + \frac{N_t(j-1)}{J}\right)\Delta t$\;
                \If {$i < N$} {
                  Guess $w_i^{[0]}(t_\ell)$\;
                }
                \If {$i > 1 $} {
                  Guess $w_{i-1}^{[0]}(t_\ell)$\;
                }
              }

            }
            \For{$\ell=1$ \KwTo $N_t/J$}{
              Solve equation~\eqref{eqn:nnwr_dirichlet} for $u_i^{[k]}(x,t_\ell)$\;
              Compute jump in Neumann data along interfaces\;
            }
            Send Neumann data\;
          }
          Check for convergence.  If converged, {\tt break}\;
        } \Else {
          \tcc{Auxiliary update}
          \For{$j=1$ \KwTo $J$}{
            Receive Neumann data\;
            \For{$\ell = 1$ \KwTo $N_t/J$}{
              $t_\ell = \left(\ell + \frac{N_t(j-1)}{J}\right)\Delta t$\;
              Solve equation~\eqref{eqn:nnwr_neumann} for $\psi_i^{[k]}(x,t_\ell)$\;
              Solve for $w_i^{[k]}(t_\ell)$ using equation~\eqref{eqn:nnwr_dirichlet_update}\;             
            }
            \If{$k<K$} {
              Send Dirichlet data\;
            }
          }
          Check for convergence.  If converged, {\tt break}\;
        }
      }
    }
  }
  \caption{This pipeline implementation of the NNWR algorithm is able
    to utilize $2NK$ computing cores if the domain is broken into $N$
    non-overlapping subdomains and $K$ full iterates are used,
    provided $J\ge2K$.}
  \label{alg:pipeline_NNWR}  
\end{algorithm}

\begin{remark}
  Unlike the classical NNWR implementation where one iterates until
  convergence, Algorithm~\ref{alg:pipeline_NNWR} requires a specification of
  $K$, the number of waveform iterates to be computed.  One can use a
  priori error estimates to pick $K$ intelligently, but it is likely
  that there will either be wasted work (if convergence is obtained
  for $k<K$) or an unconverged solution, if $K$ is not large enough.
  All is not lost however if $K$ is not large enough, since the
  computation can be restarted using the most accurate Dirichlet
  traces.  Also, in support of dynamic resource allocation and fault
  resiliency, there is active development to allow MPI communicators
  to expand or shrink.  However, there are no plans for adopting
  dynamic MPI communicators into MPI standards in the near future.
  \label{remark:pipeline_nnwr}
\end{remark}

\begin{remark}
  The pipeline NNWR implementation has a start-up and shut-down phase
  before multiple waveform iterates can be computed in parallel, i.e.,
  at the start and end of the computation, some processing cores sit
  idle.  For $K$ full waveform iterates and $J$ time blocks, the
  peak theoretical parallel efficiency is
  \begin{align}
    \begin{cases}
      \frac{2K}{2K+J-1}, & {\rm if}\quad 2K \ge J,\\
      \frac{J}{2K+J-1}, & {\rm if}\quad 2K<J.
    \end{cases}
    \label{eqn:nnwr_theoretical_efficiency}
  \end{align}
\end{remark}

\subsection{Numerical Experiments}
\label{sec:NNWR_numerics}
As mentioned in the introduction, the discretization and convergence
of NNWR algorithms have already been discussed and established for
parabolic and hyperbolic PDEs \cite{GKM3,GKM1,BankThes}.  This section
is concerned with the efficacy of the pipeline implementation, where
many waveform iterates can be concurrently computed.  The heat
equation is solved, $u_t = u_{xx}$, subject to the initial conditions
$u(x,0)=(x-0.5)^2-0.25$ and homogeneous Dirichlet boundary data.  The
computational domain, $\Omega\times[0,T]$, is $[0,1]\times [0,0.1]$.
A centered finite difference is used to approximate the spatial
operator, and a backward Euler integration is used.  Each subdomain
solves the linear system for each time advance by performing a
backwards and forwards substitution using a pre-factored LU
decomposition of the corresponding matrix.  The {\tt C} code,
available on the author's website, uses MPI parallelism.  The reported
numerical experiments were performed using the Stampede supercomputer
at the Texas Advanced Computing Center.  The resources were provided
through the NSF-supported Extreme Science and Engineering Discovery
Environment (XSEDE) program.

\begin{figure}
  \centering
  \begin{tikzpicture}
    \begin{axis}[
      xlabel={$J$: Number of time blocks},
      xmode = log,
      xmin=5,xmax=8000,
      xtick={8,32,160,800,4000},
      xticklabels={8,32,160,800,4000},
      height=0.4\textwidth,
      legend pos = outer north east,
      ylabel={Efficiency},
      legend style={cells={align=left}}]
      
      \addplot coordinates {
        (8   , 0.43)
        (16  , 0.52)
        (32  , 0.59)
        (64, 0.90)
        (128, 0.95)
        (256, 0.97)
        (512, 0.98)
        (1024, 0.99)
        (2048, 0.99)
        (4096, 1.00)
        (8192, 0.99)
      };
      \addlegendentry{Actual}
      
      \addplot[mark=none,dashed] coordinates {
        (8   , 0.53)
        (16  , 0.70)
        (32  , 0.82)
        (64  , 0.90)
        (128, 0.95)
        (256, 0.98)
        (512, 0.99)
        (1024, 0.99)
        (2048, 1.00)
        (4096, 1.00)
        (8192, 1.00)
      };
      \addlegendentry{Theoretical Peak\\(no communication)}
      
    \end{axis}
  \end{tikzpicture}
  \caption{Efficiency of the pipeline NNWR implementation as a
    function of the number of time blocks, $J$.  For a large number of
    time blocks, $J$, the algorithm is able to utilize all processors
    in a pipeline fashion for a larger percentage of the computation,
    leading to higher efficiency.}
  \label{fig:nnwr_vary_J}
\end{figure}
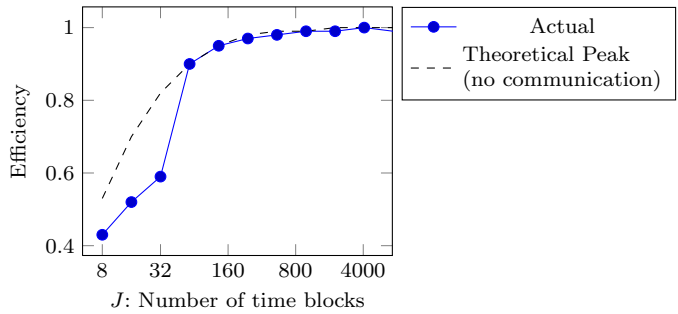

The first numerical experiment validates the qualitative behavior of
the theoretical peak efficiency,
equation~\eqref{eqn:nnwr_theoretical_efficiency}, as the number of
time blocks, $J$, is varied.  The number of full waveform iterates is
fixed at $K=4$, the number of subdomains is fixed at $N=8$, the
spatial and temporal discretizations are fixed at $N_x = 32000$,
$N_t=8192$.  A total of $64$ processing cores are used in this first
experiment.  Figure~\ref{fig:nnwr_vary_J} displays the expected
behavior -- for a small number of time blocks, processors sit idle for
a larger percentage of time, leading to poor efficiency.  For large
number of time blocks, the pipeline implementation attains close to
theoretical peak efficiency, indicating that the communication
overhead is negligible.  Here, an efficiency close to 1 means that
the pipeline NNWR implementation with $2KN$ processing cores is able
to compute $K$ full waveform iterates $2K$ times faster than the the
classical NNWR implementation using $N$ processors.  The data used to
generate Figure~\ref{fig:nnwr_vary_J} is summarized in Appendix,
Table~\ref{tbl:efficiency_nn}.

\begin{figure}
  \centering
  \begin{tikzpicture}
    \begin{semilogxaxis}[xlabel={Waveform iterates $K$},
        xmin=0,xmax=35,
        xtick={1,2,4,8,16,32,64},
        legend pos = outer north east,
        ylabel={Walltime (s)},
        log basis x={2},
        height=0.4\textwidth
      ]
      
      \addplot coordinates {
        (1, 211)
        (2, 209)
        (4, 214)
        (8, 212)
        (16, 212)
        (32, 212)
        (64, 212)
      };
    \end{semilogxaxis}
  \end{tikzpicture}
  \caption{Weak Scaling for Pipeline NNWR: Wall time vs Waveform
    iterations.  The pipeline implementation scales weakly with almost
    no overhead, i.e., with $2NK$ processing cores, we can compute $K$
    iterations in almost the same walltime as computing one iteration
    with $2N$ processing cores.}
  \label{fig:nnwr_weak_scaling}
\end{figure}
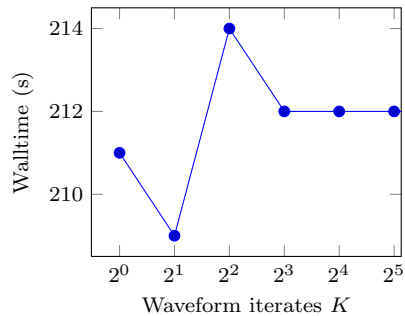

A weak scaling study is performed for the second numerical experiment.
For a fixed spatial, temporal and time-block discretization, the
number of processors is varied to compute a differing number of
waveform iterates.  In the experiment, we pick $N=8$ subdomains with
$N_x=32000$ and $N_t = 8192$, $J=1024$ time blocks and $2KN$ processor
cores.  Figure~\ref{fig:nnwr_weak_scaling} shows that the
pipeline NNWR implementation is able to scale weakly with very good 
efficiency.
The data used to generate Figure~\ref{fig:nnwr_weak_scaling} is
summarized in the Appendix, Table~\ref{tbl:weak_efficiency_nn}.

\section{Dirichlet--Neumann Waveform Relaxation (DNWR)}
\label{sec:DNWR}

The DNWR method \cite{GKM1,GKM2} can be implemented in various
arrangements, in terms of how the Dirichlet and Neumann transmission
conditions are enforced along artificial boundaries.  Pipeline
parallelism can be applied to all the variants, although the efficacy
will vary with each variant.  This section discusses the pipeline
implementation applied to the arrangement proposed by Funaro,
Quateroni and Zanolli \cite{MR972451}.  Using the notation for the
subdomain discretization in Section~\ref{sec:NNWR}, the DNWR algorithm
also requires an initial guess for the solution on the subdomain
boundaries, $w_{i}^{[0]}(t), i=1,\ldots N$, as well as a selection of
the initial subdomain, $\Omega_m$, with $1\le m \le N$; picking $m$ to
be the middle subdomain, e.g., $m=\lceil{N/2\rceil}$ will reduce the
start-up overhead for the DNWR algorithm.  Each space-time subdomain
solves one of the following subproblems.  In subdomain $m$, a
Dirichlet subproblem is solved,
\begin{subequations}
  \label{eqn:dnwr_dd}
  \begin{align}
    \partial_{t}u_{m}^{[k]} - \mathcal{L}u_{m}^{[k]}  =  f,
    &\quad (x,t) \in \Omega_m,\\
    u_{m}^{[k]}(x,0)  =  u_{0}(x),
    &\quad x \in (x_{m-1},x_m),\\
    u_{m}^{[k]}(x_{m-1},t) =
    &\begin{dcases}
       g_l(t) & \text{if } m = 1 \\
       w_{m-1}^{[k-1]}(t) & \text{otherwise}
    \end{dcases}, \\
    u_{m}^{[k]}(x_{m},t) =
    &\begin{dcases}
      g_r(t) & \text{if } m = N \\
      w_{m}^{[k-1]}(t) & \text{otherwise}
    \end{dcases}, 
  \end{align}
\end{subequations}
For subdomains to the left of subdomain $m$, i.e., $i<m$, a
Dirichlet-Neumann subproblem is solved,
\begin{subequations}
  \label{eqn:dnwr_dn}
  \begin{align}
    \partial_{t}u_{i}^{[k]} - \mathcal{L}u_{i}^{[k]}  =  f,
    &\quad (x,t) \in \Omega_i,\\
    u_{i}^{[k]}(x,0)  =  u_{0}(x),
    &\quad x \in (x_{i-1},x_i),\\
    u_{i}^{[k]}(x_{i-1},t) =
    &\begin{dcases}
       g_l(t) & \text{if } i = 1 \\
       w_{i-1}^{[k-1]}(t) & \text{otherwise}
    \end{dcases}, \\
    \partial_x u_{i}^{[k]}(x_{i},t) &=
    \partial_x u_{i+1}^{[k]}(x_{i},t) 
  \end{align}
\end{subequations}
For subdomains to the right of subdomain $m$, i.e., $i > m$, a
Neumann--Dirichlet subproblem is solved,
\begin{subequations}
  \label{eqn:dnwr_nd}
  \begin{align}
    \partial_{t}u_{i}^{[k]} - \mathcal{L}u_{i}^{[k]}  =  f,
    &\quad (x,t) \in \Omega_i,\\
    u_{i}^{[k]}(x,0)  =  u_{0}(x),
    &\quad x \in (x_{i-1},x_i),\\
    \partial_x u_{i}^{[k]}(x_{i-1},t) &=
    \partial_x u_{i-1}^{[k]}(x_{i-1},t) \\
     u_{i}^{[k]}(x_{i},t) =
    &\begin{dcases}
       g_r(t) & \text{if } i = N \\
       w_{i}^{[k-1]}(t) & \text{otherwise}
    \end{dcases}, 
  \end{align}
\end{subequations}
After the subdomain problems are solved, the Dirichlet traces are updated,
\begin{subequations}
  \label{eqn:dnwr_dirichlet_update}
  \begin{align}
    w_{i}^{[k]}(t) &= \theta u_{i}^{[k]}(x_i,t) + (1-\theta)w_{i}^{[k-1]}(t), \quad i < m,\\
    w_{i}^{[k]}(t) &= \theta u_{i+1}^{[k]}(x_i,t) + (1-\theta)w_{i}^{[k-1]}(t), \quad i \ge m.
  \end{align}
\end{subequations}
\begin{remark}
  The DNWR algorithm converges and  the rates of convergence have been
  analyzed for $N=2$ in \cite{GKM1} and for $N>2$ in \cite{GKM2}.  For the
  case $N>2$, $\theta=1/2$ and $m=\lceil N/2 \rceil$, the DNWR
  decomposition of the heat equation in $\mathbb{R}^1$ results in the
  following superlinear convergence estimate:
  \begin{align*}
    \max_{1\leq i\leq N-1}\| w_{i}^{[k]}\|_{L^{\infty}(0,T)}
    \leq\left(N-4+\frac{2h_{\max}}{h_m}\right)^{k}
    \erfc \left(\frac{kh_{\min}}{2\sqrt{T}}\right)
    \max_{1\leq i\leq N-1} \| w_{i}^{[0]}\|_{L^{\infty}(0,T)},
  \end{align*}
  where $h_i=x_i-x_{i-1}$, $\displaystyle h_{\max}:=\max_{1\leq i\leq N}h_i$ and
  $\displaystyle h_{\min}:=\min_{1\leq i\leq N}h_i$.
\end{remark}

\subsection{Classical DNWR implementation}

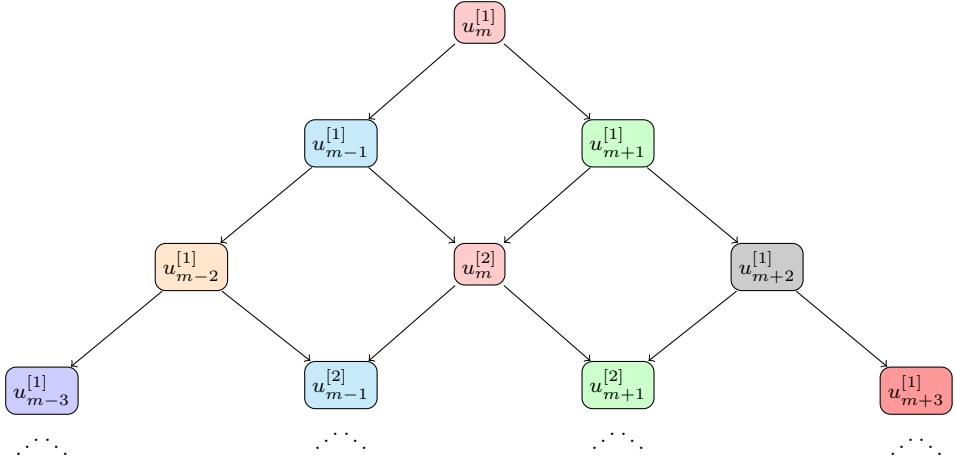
\begin{figure}
  \begin{tikzpicture}

  \tikzstyle{p1} = [rectangle, draw, fill=blue!20, rounded corners]
  \tikzstyle{p2} = [rectangle, draw, fill=orange!20, rounded corners]
  \tikzstyle{p3} = [rectangle, draw, fill=cyan!20, rounded corners]
  \tikzstyle{p4} = [rectangle, draw, fill=red!20, rounded corners]
  \tikzstyle{p5} = [rectangle, draw, fill=green!20, rounded corners]
  \tikzstyle{p6} = [rectangle, draw, fill=black!20, rounded corners]
  \tikzstyle{p7} = [rectangle, draw, fill=red!40, rounded corners]

  \node[p4] (um1) {$u_m^{[1]}$};

  \node[p3,below left= of um1] (ul1) {$u_{m-1}^{[1]}$};
  \node[p5,below right= of um1] (un1) {$u_{m+1}^{[1]}$};

  \node[p2,below left= of ul1] (uk1) {$u_{m-2}^{[1]}$};
  \node[p4,below right= of ul1] (um2) {$u_{m}^{[2]}$};
  \node[p6,below right= of un1] (uo1) {$u_{m+2}^{[1]}$};

  \node[p1,below left= of uk1] (uj1) {$u_{m-3}^{[1]}$};
  \node[p3,below left= of um2] (ul2) {$u_{m-1}^{[2]}$};
  \node[p5,below right= of um2] (un2) {$u_{m+1}^{[2]}$};
  \node[p7,below right= of uo1] (up1) {$u_{m+3}^{[1]}$};

  \draw[->] (um1)--(ul1);
  \draw[->] (um1)--(un1);

  \draw[->] (ul1)--(uk1);
  \draw[->] (ul1)--(um2);
  \draw[->] (un1)--(um2);
  \draw[->] (un1)--(uo1);

  \draw[->] (uk1)--(uj1);
  \draw[->] (uk1)--(ul2);
  \draw[->] (um2)--(ul2);
  \draw[->] (um2)--(un2);
  \draw[->] (uo1)--(up1);
  \draw[->] (uo1)--(un2);

  \node[below = of un2,yshift = 1cm] {$\iddots\ddots$};
  \node[below = of up1,yshift = 1cm] {$\iddots\ddots$};
  \node[below = of uj1,yshift = 1cm] {$\iddots\ddots$};
  \node[below = of ul2,yshift = 1cm] {$\iddots\ddots$};
\end{tikzpicture}
  \caption{DNWR implemented classically. The different colors
    represent $N$ different processing cores (inefficiently) used to
    compute the waveform iterates in different space--time subdomains.
  }
  \label{fig:classical_dnwr}
\end{figure}

Equations \eqref{eqn:dnwr_dd} -- \eqref{eqn:dnwr_dn} --
\eqref{eqn:dnwr_nd} -- \eqref{eqn:dnwr_dirichlet_update} can be
``decoupled'' by solving the space--time subproblems in a specific
order.  Recall the decomposition of the space-time domain as described
in Figure~\ref{fig:subdomains}.  One starts by computing
$u_m^{[1]}(x,t)$ in $\Omega_m$ and transmitting the computed boundary
conditions before computing $u_{m-1}^{[1]}(x,t)$ and
$u_{m+1}^{[1]}(x,t)$ in the neighboring subdomains.  The updated
boundary conditions are transmitted, and then $u_{m-2}^{[1]}(x,t)$,
$u_m^{[2]}(x,t)$ and $u_{m+1}^{[1]}(x,t)$ are simultaneously computed.
If optimal parallel efficiency is not a concern, $N$ processes can be
spawned as shown in Algorithm~\ref{alg:DNWR_naive}.  Not all processes
can be utilized simultaneously however, because transmission
conditions needed in lines 6 through 19 may not be available until
other processes have computed and transmitted the required
information.  Although not a practical algorithm, this pseudo code is
instructive because it outlines the flow of information in an easy to
read fashion. Figure~\ref{fig:classical_dnwr} shows the flow of
information for the first few steps of the DNWR algorithm, implemented
classically using $N$ processors and $N$ subdomains.

\scalebox{0.95}{
\begin{algorithm}[H]
  \ForPar{$i=1$ \KwTo $N$}{
    $m = nprocs = \lceil N/2 \rceil$\;
    \If {$ m \le i < N $} {
      Guess $w_i^{[0]}(t_\ell), \quad \ell = 1,\ldots N_t$ \;
    }
    \If {$1 < i \le m$} {
      Guess $w_{i-1}^{[0]}(t_\ell), \quad \ell = 1,\ldots N_t$ \;
    }
    \For{$k=1$ \KwTo $K$}{
      \tcc{receive boundary information}
      \Switch{i}{
        \Case{$i=m$} {
          \If {$k > 1$} {
            Receive Dirichlet data from neighbors (if they exist)\;
            Update Dirichlet traces using equation~\eqref{eqn:dnwr_dirichlet_update}\;
          }
        }
        \Case {$i<m$} {
          Receive Neumann data from right neighbor\;
          \If {$k > 1$ \bf{and} $i>1$} {
            Receive Dirichlet data from left neighbors\;
            Update Dirichlet traces using equation~\eqref{eqn:dnwr_dirichlet_update}\;
          }

        }
        \Case {$i>m$} {
          Receive Neumann data from left neighbor\;
          \If {$k > 1$ \bf{and} $i<N$} {
            Receive Dirichlet data from right neighbors\;
            Update Dirichlet traces using equation~\eqref{eqn:dnwr_dirichlet_update}\;
          }
        }
      }

      \tcc{Solve space time problem}
      \For {$\ell = 1$ \KwTo $N_t$} {
        \Switch{i}{
          \Case{$i=m$} {
            Solve equation~\eqref{eqn:dnwr_dd} for $u_i^{[k]}(x,t_\ell)$\;
          }
          \Case {$i<m$} {
            Solve equation~\eqref{eqn:dnwr_nd} for $u_i^{[k]}(x,t_\ell)$\;
          }
          \Case {$i>m$} {
            Solve equation~\eqref{eqn:dnwr_dn} for $u_i^{[k]}(x,t_\ell)$\;
          }
        }
      } 

      \tcc{Send boundary conditions}
      \Switch{i}{
        \Case{$i=m$} {
          Send Neumann data to neighbors (if they exist)\;
        }
        \Case {$i<m$} {
          \If {$i>1$} {
            Send Neumann/convergence data to left neighbor\;
          }
          \If {$k < K$} {
            Send Dirichlet data to right neighbor\;
          }
        }
        \Case {$i>m$} {
          \If {$i<N$} {
            Send Neumann/convergence data to right neighbor\;
          }
          \If {$k < K$} {
            Send Dirichlet data to left neighbor\;
          }
        }
      }
      Check for convergence.  If converged, {\tt break}\;
    }
  }
  
  \caption{A classical (naive) implementation of the DNWR algorithm
    using $N$ processes, where $K$ is the number of waveform iterates
    and $N$ is the number of non-overlapping spatial subdomains.}
  \label{alg:DNWR_naive}  
\end{algorithm}
}

\newpage
\begin{remark}
  By sending appropriate convergence flags in lines 37 and 42 of
  Algorithm~\ref{alg:DNWR_naive}, the processors computing the
  solution to the space time block in $\Omega_1$ and $\Omega_N$ are
  able to determine if the algorithm has converged.
\end{remark}

To construct practical DNWR algorithms using a classical
implementation, recall from the previous section that selecting
$m=\lceil N/2\rceil$ reduces the startup overhead for the DNWR
algorithm, thereby improving parallel efficacy.  For the remainder of
this paper, assume that $m=\lceil N/2\rceil$.  The classical DNWR
algorithm depends on the ratio between the number of subdomains and
the number of waveform iterates being computed.  If $2K \ge \lceil N/2
\rceil$, the DNWR algorithm is able to utilize $\lceil N/2 \rceil$
processing cores -- each processor computes the solution for two
spatial subdomains.
Algorithm~\ref{alg:DNWR} outlines a pseudo code of the DNWR
implemented classically.  Algorithm~\ref{alg:DNWR_case2} in the
Appendix discusses the implementation for the case $2K < \lceil N/2
\rceil$.  For both cases, a total of $(N-1)(2K-1)$ messages are
needed, with each message containing $N_t$ words, where $N_t$ is the
number of time steps.

\begin{remark}
   Unlike the classical implementation of the NNWR algorithm, the
   classical implementation of the DNWR algorithm requires
   specification of $K$, the number of waveform iterates desired.
   This is because each spatial subdomain might compute different
   waveform iterates simultaneously.  Similar to
   Remark~\ref{remark:pipeline_nnwr}, a priori specification of $K$
   might result in an unconverged solution if the prescribed $K$ is
   not large enough, or wasted work if convergence is obtained for
   $k<K$ iterations.
\end{remark}

\begin{algorithm}
  \KwIn{$N$: \# subdomains; $K$: \# waveform iterates}
  $m = \lceil N/2 \rceil$.\;
  
  \ForPar{$p=1$ \KwTo $nprocs$}{
    Set $i = 2(p-1)$\;
    \If {$i>0$ \&\&  $i\le m$} {
      Guess $w_{i}^{[0]}(t_\ell),\quad \ell =1\ldots,N_t$ \;
    }
    Set $i = 2p$\;
    \If {$i<N$ \&\&  $i\ge m$} {
      Guess $w_{i}^{[0]}(t_\ell),\quad \ell =1\ldots,N_t$ \;
    }
    Set $i = 2p-1$\;
    \If {$i<N$} {
      Guess $w_{i}^{[0]}(t_\ell),\quad \ell =1\ldots,N_t$ \;
    }

    \For{$k=1$ \KwTo $K$}{
      \Switch{$p$}{
        \Case{$2p < m$} {
          Set $i = 2p$\;
        }
        \Case{$(2p-1) > m$} {
          Set $i = 2p-1$\;
        }
        \Case{$2p = m$}{
          Set $i = 2p$\;
        }
        \Other {
          Set $i = 2p-1$\;
        }
      }
      \If {$i\in\{1,\ldots,N\}$} {
        \tcc{receive boundary info: lines 8--22 in Alg.~\ref{alg:DNWR_naive}}
        \tcc{Solve space time pde: lines 23--30 in Alg.~\ref{alg:DNWR_naive}}
        \tcc{Send boundary info: lines 31--43 in Alg.~\ref{alg:DNWR_naive}}
      }

      \Switch{$p$}{
        \Case{$2p < m$} {
          Set $i = 2p-1$\;
        }
        \Case{$(2p-1) > m$} {
          Set $i = 2p$\;
        }
        \Case{$2p = m$}{
          Set $i = 2p-1$\;
        }
        \Other {
          Set $i = 2p$\;
        }
      }
      \If {$i\in\{1,\ldots,N\}$} {
        \tcc{receive boundary info: lines 8--22 in Alg.~\ref{alg:DNWR_naive}}
        \tcc{Solve space time pde: lines 23--30 in Alg.~\ref{alg:DNWR_naive}}
        \tcc{Send boundary info: lines 31--43 in Alg.~\ref{alg:DNWR_naive}}
      }
      Check for convergence.  If converged, {\tt break}\;
    }
  }
  \caption{A classical implementation of the DNWR algorithm using
    $\lceil N/2 \rceil$ processes, where $N$ is the number of
    non-overlapping spatial subdomains.}
  \label{alg:DNWR}  
\end{algorithm}

\subsection{Pipeline DNWR implementation}
\label{sec:pDNWR}

Similar to the pipeline NNWR implementation in
Section~\ref{sec:pipeline_NNWR}, a pipeline implementation of the DNWR
algorithm allows for higher concurrency at the expense of increased
communication: multiple waveform iterates are simultaneously
evaluated (in $\cup_j \Omega_{ij}$) after an initial startup cost.  As
before, both the classical implementation and the pipeline
implementation of the DNWR algorithm will result in the same numerical
solution.  Figure~\ref{fig:pipeline_dnwr} shows a graphical
representation of a five-domain, two waveform iterate example.  For a
general decomposition as shown in Figure~\ref{fig:pipeline_domains},
the pipeline DNWR implementation is outlined in
Algorithm~\ref{alg:pipeline_DNWR}. If the number of time blocks, $J$,
are sufficiently large, the pipeline DNWR algorithm is able to utilize
$NK$ processing cores.  The pipeline DNWR implementation requires a
total of $J(N-1)(2K-1)$ messages, with each message containing $N_t/J$
words.  Since the numerical results in Section~\ref{sec:NNWR_numerics}
indicate that the communication overhead of utilizing many time
blocks $J$ is negligible, we restrict our discussion to the case where
$J$ is sufficiently large, $J > \lceil N/2 \rceil + (2K - 1)$.

\begin{landscape}
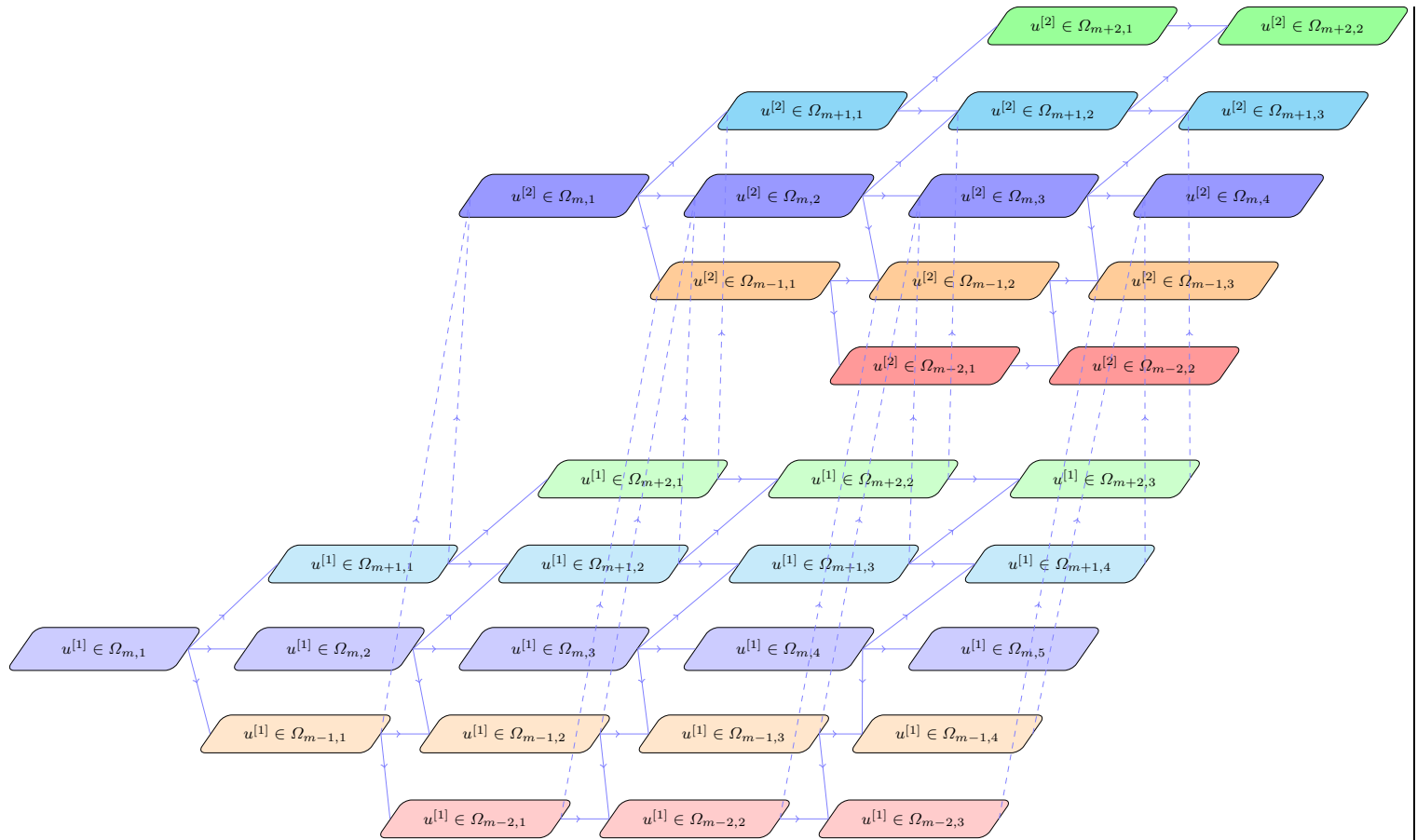
\begin{figure}
  \centering
  \scalebox{0.8}{
    \begin{tikzpicture}

  trapezium, draw, minimum width=3cm,
  trapezium left angle=120, trapezium right angle=60
  
  \tikzstyle{p1} = [trapezium, trapezium left angle=55, trapezium right angle=125, draw,
    fill=blue!20, rounded corners, minimum width=3.5cm]
  \tikzstyle{p2} = [trapezium, trapezium left angle=55, trapezium right angle=125, draw,
    fill=orange!20, rounded corners, minimum width=3.5cm]
  \tikzstyle{p3} = [trapezium, trapezium left angle=55, trapezium right angle=125, draw,
    fill=cyan!20, rounded corners, minimum width=3.5cm]
  \tikzstyle{p4} = [trapezium, trapezium left angle=55, trapezium right angle=125, draw,
    fill=red!20, rounded corners, minimum width=3.5cm]
  \tikzstyle{p5} = [trapezium, trapezium left angle=55, trapezium right angle=125, draw,
    fill=green!20, rounded corners, minimum width=3.5cm]

  \tikzstyle{p6} = [trapezium, trapezium left angle=55, trapezium right angle=125, draw,
    fill=blue!40, rounded corners, minimum width=3.5cm]
  \tikzstyle{p7} = [trapezium, trapezium left angle=55, trapezium right angle=125, draw,
    fill=orange!40, rounded corners, minimum width=3.5cm]
  \tikzstyle{p8} = [trapezium, trapezium left angle=55, trapezium right angle=125, draw,
    fill=cyan!40, rounded corners, minimum width=3.5cm]
  \tikzstyle{p9} = [trapezium, trapezium left angle=55, trapezium right angle=125, draw,
    fill=red!40, rounded corners, minimum width=3.5cm]
  \tikzstyle{p10} = [trapezium, trapezium left angle=55, trapezium right angle=125, draw,
    fill=green!40, rounded corners, minimum width=3.5cm]

  \node[p1] at (0,0) (u1m1) {$u^{[1]} \in \Omega_{m,1}$};
  \node[p1] at (4,0) (u1m2) {$u^{[1]} \in \Omega_{m,2}$};
  \node[p1] at (8,0) (u1m3) {$u^{[1]} \in \Omega_{m,3}$};
  \node[p1] at (12,0) (u1m4) {$u^{[1]} \in \Omega_{m,4}$};
  \node[p1] at (16,0) (u1m5) {$u^{[1]} \in \Omega_{m,5}$};

  \node[p2] at (3.4, -1.5) (u1l1) {$u^{[1]} \in \Omega_{m-1,1}$};
  \node[p2] at (7.3, -1.5) (u1l2) {$u^{[1]} \in \Omega_{m-1,2}$};
  \node[p2] at (11.2, -1.5) (u1l3) {$u^{[1]} \in \Omega_{m-1,3}$};
  \node[p2] at (15, -1.5) (u1l4) {$u^{[1]} \in \Omega_{m-1,4}$};

  \node[p3] at (4.6, 1.5) (u1n1) {$u^{[1]} \in \Omega_{m+1,1}$};
  \node[p3] at (8.7, 1.5) (u1n2) {$u^{[1]} \in \Omega_{m+1,2}$};
  \node[p3] at (12.8, 1.5) (u1n3) {$u^{[1]} \in \Omega_{m+1,3}$};
  \node[p3] at (17, 1.5) (u1n4) {$u^{[1]} \in \Omega_{m+1,4}$};

  \node[p4] at (6.6, -3) (u1k1) {$u^{[1]} \in \Omega_{m-2,1}$};
  \node[p4] at (10.5, -3) (u1k2) {$u^{[1]} \in \Omega_{m-2,2}$};
  \node[p4] at (14.4, -3) (u1k3) {$u^{[1]} \in \Omega_{m-2,3}$};

  \node[p5] at (9.4, 3) (u1o1) {$u^{[1]} \in \Omega_{m+2,1}$};
  \node[p5] at (13.5, 3) (u1o2) {$u^{[1]} \in \Omega_{m+2,2}$};
  \node[p5] at (17.8, 3) (u1o3) {$u^{[1]} \in \Omega_{m+2,3}$};

  \node[p6] at (8,8) (u2m1) {$u^{[2]} \in \Omega_{m,1}$};
  \node[p6] at (12,8) (u2m2) {$u^{[2]} \in \Omega_{m,2}$};
  \node[p6] at (16,8) (u2m3) {$u^{[2]} \in \Omega_{m,3}$};
  \node[p6] at (20,8) (u2m4) {$u^{[2]} \in \Omega_{m,4}$};

  \node[p7] at (11.4, 6.5) (u2l1) {$u^{[2]} \in \Omega_{m-1,1}$};
  \node[p7] at (15.3, 6.5) (u2l2) {$u^{[2]} \in \Omega_{m-1,2}$};
  \node[p7] at (19.2, 6.5) (u2l3) {$u^{[2]} \in \Omega_{m-1,3}$};

  \node[p8] at (12.6, 9.5) (u2n1) {$u^{[2]} \in \Omega_{m+1,1}$};
  \node[p8] at (16.7, 9.5) (u2n2) {$u^{[2]} \in \Omega_{m+1,2}$};
  \node[p8] at (20.8, 9.5) (u2n3) {$u^{[2]} \in \Omega_{m+1,3}$};

  \node[p9] at (14.6, 5) (u2k1) {$u^{[2]} \in \Omega_{m-2,1}$};
  \node[p9] at (18.5, 5) (u2k2) {$u^{[2]} \in \Omega_{m-2,2}$};

  \node[p10] at (17.4, 11) (u2o1) {$u^{[2]} \in \Omega_{m+2,1}$};
  \node[p10] at (21.5, 11) (u2o2) {$u^{[2]} \in \Omega_{m+2,2}$};

  
  \begin{scope}[color=blue!50,decoration={
        markings,
        mark=at position 0.4 with {\arrow{>}}}
    ]
    
    \draw[postaction={decorate}] (u1m1.east) -- (u1m2.west);
    \draw[postaction={decorate}] (u1m1.east) -- (u1n1.west);
    \draw[postaction={decorate}] (u1m1.east) -- (u1l1.west);
    
    \draw[postaction={decorate}] (u1m2.east) -- (u1m3.west);
    \draw[postaction={decorate}] (u1m2.east) -- (u1n2.west);
    \draw[postaction={decorate}] (u1m2.east) -- (u1l2.west);
    \draw[postaction={decorate}] (u1n1.east) -- (u1n2.west);
    \draw[postaction={decorate}] (u1n1.east) -- (u1o1.west);
    \draw[postaction={decorate}] (u1l1.east) -- (u1l2.west);
    \draw[postaction={decorate}] (u1l1.east) -- (u1k1.west);

    \draw[postaction={decorate}] (u1m3.east) -- (u1m4.west);
    \draw[postaction={decorate}] (u1m3.east) -- (u1n3.west);
    \draw[postaction={decorate}] (u1m3.east) -- (u1l3.west);
    \draw[postaction={decorate}] (u1n2.east) -- (u1n3.west);
    \draw[postaction={decorate}] (u1n2.east) -- (u1o2.west);
    \draw[postaction={decorate}] (u1l2.east) -- (u1l3.west);
    \draw[postaction={decorate}] (u1l2.east) -- (u1k2.west);

    \draw[postaction={decorate}] (u1o1.east) -- (u1o2.west);
    \draw[postaction={decorate}] (u1o2.east) -- (u1o3.west);

    \draw[postaction={decorate}] (u1k1.east) -- (u1k2.west);
    \draw[postaction={decorate}] (u1k2.east) -- (u1k3.west);

    \draw[postaction={decorate}] (u1m4.east) -- (u1m5.west);
    \draw[postaction={decorate}] (u1m4.east) -- (u1n4.west);
    \draw[postaction={decorate}] (u1m4.east) -- (u1l4.west);
    \draw[postaction={decorate}] (u1n3.east) -- (u1n4.west);
    \draw[postaction={decorate}] (u1n3.east) -- (u1o3.west);
    \draw[postaction={decorate}] (u1l3.east) -- (u1l4.west);
    \draw[postaction={decorate}] (u1l3.east) -- (u1k3.west);

    \draw[postaction={decorate}] (u2m1.east) -- (u2m2.west);
    \draw[postaction={decorate}] (u2m1.east) -- (u2n1.west);
    \draw[postaction={decorate}] (u2m1.east) -- (u2l1.west);
    
    \draw[postaction={decorate}] (u2m2.east) -- (u2m3.west);
    \draw[postaction={decorate}] (u2m2.east) -- (u2n2.west);
    \draw[postaction={decorate}] (u2m2.east) -- (u2l2.west);
    \draw[postaction={decorate}] (u2n1.east) -- (u2n2.west);
    \draw[postaction={decorate}] (u2n1.east) -- (u2o1.west);
    \draw[postaction={decorate}] (u2l1.east) -- (u2l2.west);
    \draw[postaction={decorate}] (u2l1.east) -- (u2k1.west);

    \draw[postaction={decorate}] (u2m3.east) -- (u2m4.west);
    \draw[postaction={decorate}] (u2m3.east) -- (u2n3.west);
    \draw[postaction={decorate}] (u2m3.east) -- (u2l3.west);
    \draw[postaction={decorate}] (u2n2.east) -- (u2n3.west);
    \draw[postaction={decorate}] (u2n2.east) -- (u2o2.west);
    \draw[postaction={decorate}] (u2l2.east) -- (u2l3.west);
    \draw[postaction={decorate}] (u2l2.east) -- (u2k2.west);

    \draw[postaction={decorate}] (u2o1.east) -- (u2o2.west);
    \draw[postaction={decorate}] (u2k1.east) -- (u2k2.west);
    
    \draw[postaction={decorate},dashed] (u1n1.east) -- (u2m1.west);
    \draw[postaction={decorate},dashed] (u1l1.east) -- (u2m1.west);
    \draw[postaction={decorate},dashed] (u1n2.east) -- (u2m2.west);
    \draw[postaction={decorate},dashed] (u1l2.east) -- (u2m2.west);
    \draw[postaction={decorate},dashed] (u1n3.east) -- (u2m3.west);
    \draw[postaction={decorate},dashed] (u1l3.east) -- (u2m3.west);
    \draw[postaction={decorate},dashed] (u1n4.east) -- (u2m4.west);
    \draw[postaction={decorate},dashed] (u1l4.east) -- (u2m4.west);

    \draw[postaction={decorate},dashed] (u1k1.east) -- (u2l1.west);
    \draw[postaction={decorate},dashed] (u1k2.east) -- (u2l2.west);
    \draw[postaction={decorate},dashed] (u1k3.east) -- (u2l3.west);

    \draw[postaction={decorate},dashed] (u1o1.east) -- (u2n1.west);
    \draw[postaction={decorate},dashed] (u1o2.east) -- (u2n2.west);
    \draw[postaction={decorate},dashed] (u1o3.east) -- (u2n3.west);

  \end{scope}
  
\end{tikzpicture}
    }
  \caption{A pipeline DNWR implementation for a five domain, two
    waveform iterate example.  Each color/shade represents a different
    processor.}
  \label{fig:pipeline_dnwr}
\end{figure}
\end{landscape}

\begin{algorithm}
  \tcc{Assumes $J > \lceil N/2 \rceil + (2K - 1)$}
  \KwIn{$N$: \# subdomains; $K$: \# waveform iterates}
  $m = \lceil N/2 \rceil$\;
  
  \ForPar{$i=1$ \KwTo $N$}{
    \ForPar{$k=1$ \KwTo $K$}{
      
      \For{$j = 1$ \KwTo $J$}{
        \tcc{Receive boundary data, lines 8--22 in Alg.~\ref{alg:DNWR_naive}}
        \For{$\ell=1$ \KwTo $N_t/J$}{
          $t_\ell = \left(\ell + \frac{N_t(j-1)}{J}\right)\Delta t$\;
          \If {$k = 1$} {
            \If {$i < N$} {
              Guess $w_i^{[0]}(t_n)$\;
            }
            \If {$i > 0 $} {
              Guess $w_{i-1}^{[0]}(t_n)$\;
            }
          }

          \tcc{Solve space-time PDE, lines 23--30, Alg.~\ref{alg:DNWR_naive}}
        }

        \tcc{Send boundary conditions, lines 31--43 in Alg.~\ref{alg:DNWR_naive}}
      }
      Check for convergence.  If converged, {\tt break}\;
    }
  }
  \caption{The pipeline implementation of the DNWR algorithm is able
    to utilize $NK$ computing cores if the domain is broken into $N$
    non-overlapping subdomains.}
  \label{alg:pipeline_DNWR}  
\end{algorithm}

\begin{remark}
  The pipeline DNWR implementation has a start-up and shut-down phase
  before multiple waveform iterates can be computed in parallel. For
  $K$ waveform iterates and $J$ time blocks, the parallel efficiency
  is
  \begin{align}
    \frac{J}{J+\lfloor{(N/2)}\rfloor+2(K-1)},
  \end{align}
  provided $J > \lceil N/2 \rceil + (2K - 1)$.
\end{remark}

\subsection{Numerical Experiment}

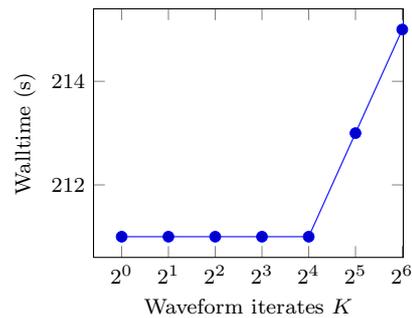
\begin{figure}
  \centering
  \begin{tikzpicture}
    \begin{semilogxaxis}[xlabel={Waveform iterates $K$},
        xmin=0,xmax=65,
        xtick={1,2,4,8,16,32,64,128},
        legend pos = outer north east,
        ylabel={Walltime (s)},
        log basis x={2},
        height=0.4\textwidth
      ]
    
      \addplot coordinates {
        (1, 211)
        (2, 211)
        (4, 211)
        (8, 211)
        (16, 211)
        (32, 213)
        (64, 215)
        (128,220)
      };
    \end{semilogxaxis}
  \end{tikzpicture}
  
  \caption{Weak Scaling for Pipeline DNWR: Wall time vs Waveform
    iterations.  The pipeline implementation scales weakly with almost
    no overhead, i.e., with $NK$ processing cores, we can compute $K$
    iterations in almost the same walltime as computing one iteration
    with $N$ processing cores.}
  \label{fig:dnwr_weak_scaling}
\end{figure}

A weak scaling study of the pipeline DNWR algorithm concludes the
discussion of implementation details related to waveform relaxation
methods. The heat equation described in
Section~\ref{sec:NNWR_numerics} is solved with a fixed spatial,
temporal and time-block discretization.  The number of processors is
varied to compute a differing number of waveform iterates.  In the
experiment, we pick $N=8$ subdomains with $N_x = 32,000$, $N_t = 8192$
and $J = 1024$ time blocks.  Figure~\ref{fig:dnwr_weak_scaling} shows
that the pipeline DNWR implementation is able to scale weakly with
very good efficiency.
The data used to generate Figure~\ref{fig:dnwr_weak_scaling} is
summarized in the Appendix, Table~\ref{tbl:weak_efficiency_dnwr}.

\section{Conclusions}
\label{sec:conclusion}

In this paper, we have reformulated the NNWR and DNWR methods to allow
for pipeline-parallel computation of the waveform iterates, after an
initial start-up cost.  The key observation is that the order of the
computations (loops) can be flipped without affecting the final
solution, at the expense of increased communication.  Theoretical
estimates for the parallel speedup and communication overhead are
presented, along with weak scaling studies to show the effectiveness
of the pipeline DNWR and NNWR algorithms.

A significant drawback, is that the pipeline implementations (and the
DNWR algorithm) require an apriori specification of the number of
waveform iterates to be computed.  Hence, it is likely that there will
either be wasted work (if convergence is obtained for fewer waveform
iterates) or an unconverged solution, if insufficient waveform
iterates were specified.  An adaptive approach is being designed and
analyzed to improve these pipeline parallel computations.

\section*{Acknowledgments}
This work utilized computational resources provided by superior, the
high-performance computing cluster @ MTU, and by the Extreme Science
and Engineering Discovery Environment (XSEDE), which is supported by
National Science Foundation grant number ACI-1053575.

\bibliography{wr}

\begin{thebibliography}{10}
\providecommand{\url}[1]{{#1}}
\providecommand{\urlprefix}{URL }
\expandafter\ifx\csname urlstyle\endcsname\relax
  \providecommand{\doi}[1]{DOI~\discretionary{}{}{}#1}\else
  \providecommand{\doi}{DOI~\discretionary{}{}{}\begingroup
  \urlstyle{rm}\Url}\fi

\bibitem{MR865945}
Bj{\o}rstad, P.E., Widlund, O.B.: Iterative methods for the solution of
  elliptic problems on regions partitioned into substructures.
\newblock SIAM J. Numer. Anal. \textbf{23}(6), 1097--1120 (1986).
\newblock \doi{10.1137/0723075}.
\newblock \urlprefix\url{http://dx.doi.org/10.1137/0723075}

\bibitem{MR992000}
Bourgat, J.F., Glowinski, R., Le~Tallec, P., Vidrascu, M.: Variational
  formulation and algorithm for trace operator in domain decomposition
  calculations.
\newblock In: Domain decomposition methods ({L}os {A}ngeles, {CA}, 1988), pp.
  3--16. SIAM, Philadelphia, PA (1989)

\bibitem{MR829613}
Bramble, J.H., Pasciak, J.E., Schatz, A.H.: An iterative method for elliptic
  problems on regions partitioned into substructures.
\newblock Math. Comp. \textbf{46}(174), 361--369 (1986).
\newblock \doi{10.2307/2007981}.
\newblock \urlprefix\url{http://dx.doi.org/10.2307/2007981}

\bibitem{MR1106455}
De~Roeck, Y.H., Le~Tallec, P.: Analysis and test of a local
  domain-decomposition preconditioner.
\newblock In: Fourth {I}nternational {S}ymposium on {D}omain {D}ecomposition
  {M}ethods for {P}artial {D}ifferential {E}quations ({M}oscow, 1990), pp.
  112--128. SIAM, Philadelphia, PA (1991)

\bibitem{MR972451}
Funaro, D., Quarteroni, A., Zanolli, P.: An iterative procedure with interface
  relaxation for domain decomposition methods.
\newblock SIAM J. Numer. Anal. \textbf{25}(6), 1213--1236 (1988).
\newblock \doi{10.1137/0725069}.
\newblock \urlprefix\url{http://dx.doi.org/10.1137/0725069}

\bibitem{MR2300292}
Gander, M.J., Halpern, L.: Optimized {S}chwarz waveform relaxation methods for
  advection reaction diffusion problems.
\newblock SIAM J. Numer. Anal. \textbf{45}(2), 666--697 (electronic) (2007).
\newblock \doi{10.1137/050642137}.
\newblock \urlprefix\url{http://dx.doi.org/10.1137/050642137}

\bibitem{MR2035001}
Gander, M.J., Halpern, L., Nataf, F.: Optimal {S}chwarz waveform relaxation for
  the one dimensional wave equation.
\newblock SIAM J. Numer. Anal. \textbf{41}(5), 1643--1681 (2003).
\newblock \doi{10.1137/S003614290139559X}.
\newblock \urlprefix\url{http://dx.doi.org/10.1137/S003614290139559X}

\bibitem{GKM1}
Gander, M.J., Kwok, F., Mandal, B.C.: Dirichlet-{N}eumann and
  {N}eumann-{N}eumann {W}aveform {R}elaxation {A}lgorithms for {P}arabolic
  {P}roblems.
\newblock Electron. Trans. Numer. Anal. \textbf{45}, 424--456 (2016)

\bibitem{GKM3}
Gander, M.J., Kwok, F., Mandal, B.C.: Dirichlet-{N}eumann and
  {N}eumann-{N}eumann {W}aveform {R}elaxation for the {W}ave {E}quation.
\newblock In: T.~Dickopf, M.J. Gander, L.~Halpern, R.~Krause, L.~Pavarino
  (eds.) Domain Decomposition Methods in Science and Engineering XXII, Lecture
  Notes in Computational Science and Engineering, pp. 501--509. Springer
  International Publishing (2016).
\newblock \doi{10.1007/978-3-319-18827-0_51}.
\newblock \urlprefix\url{http://dx.doi.org/10.1007/978-3-319-18827-0_51}

\bibitem{GKM2}
Gander, M.J., Kwok, F., Mandal, B.C.: Dirichlet-{N}eumann {W}aveform
  {R}elaxation {M}ethod for the 1{D} and 2{D} {H}eat and {W}ave {E}quations in
  {M}ultiple subdomains.
\newblock to appear  (arXiv:1507.04011)

\bibitem{MR1638096}
Gander, M.J., Stuart, A.M.: Space-time continuous analysis of waveform
  relaxation for the heat equation.
\newblock SIAM J. Sci. Comput. \textbf{19}(6), 2014--2031 (1998).
\newblock \doi{10.1137/S1064827596305337}.
\newblock \urlprefix\url{http://dx.doi.org/10.1137/S1064827596305337}

\bibitem{MR1146977}
Gear, C.W.: Waveform methods for space and time parallelism.
\newblock In: Proceedings of the {I}nternational {S}ymposium on {C}omputational
  {M}athematics ({M}atsuyama, 1990), vol.~38, pp. 137--147 (1991).
\newblock \doi{10.1016/0377-0427(91)90166-H}.
\newblock \urlprefix\url{http://dx.doi.org/10.1016/0377-0427(91)90166-H}

\bibitem{GilKel}
Giladi, E., Keller, H.: Space time domain decomposition for parabolic problems.
\newblock Tech. Rep. 97-4, Center for research on parallel computation CRPC,
  Caltech (1997)

\bibitem{Kwok}
Kwok, F.: Neumann-{N}eumann {W}aveform {R}elaxation for the {T}ime-{D}ependent
  {H}eat {E}quation.
\newblock In: J.~Erhel, M.J. Gander, L.~Halpern, G.~Pichot, T.~Sassi, O.B.
  Widlund (eds.) Domain Decomposition in Science and Engineering XXI, vol.~98,
  pp. 189--198. Springer-Verlag (2014)

\bibitem{MR1095198}
Le~Tallec, P., De~Roeck, Y.H., Vidrascu, M.: Domain decomposition methods for
  large linearly elliptic three-dimensional problems.
\newblock J. Comput. Appl. Math. \textbf{34}(1), 93--117 (1991).
\newblock \doi{10.1016/0377-0427(91)90150-I}.
\newblock \urlprefix\url{http://dx.doi.org/10.1016/0377-0427(91)90150-I}

\bibitem{LelRue}
Lelarasmee, E., Ruehli, A., Sangiovanni-Vincentelli, A.: The waveform
  relaxation method for time-domain analysis of large scale integrated
  circuits.
\newblock IEEE Trans. Compt.-Aided Design Integr. Circuits Syst. \textbf{1}(3),
  131--145 (1982)

\bibitem{Lind}
Lindel\"{o}f, E.: Sur l'application des m\'ethodes d'approximations successives
  \`a l'\'etude des int\'egrales r\'eelles des \'equations diff\'erentielles
  ordinaires.
\newblock Journal de Math\'ematiques Pures et Appliqu\'ees pp. 117--128 (1894)

\bibitem{BankThes}
Mandal, B.C.: Convergence {A}nalysis of {S}ubstructuring {W}aveform
  {R}elaxation {M}ethods for {S}pace-time {P}roblems and {T}heir {A}pplication
  to {O}ptimal {C}ontrol {P}roblems.
\newblock Ph.D. thesis, University of Geneva (2014).
\newblock \urlprefix\url{http://archive-ouverte.unige.ch/unige:46146}

\bibitem{Mandal}
Mandal, B.C.: A {T}ime-{D}ependent {D}irichlet-{N}eumann {M}ethod for the
  {H}eat {E}quation.
\newblock In: J.~Erhel, M.J. Gander, L.~Halpern, G.~Pichot, T.~Sassi, O.B.
  Widlund (eds.) Domain Decomposition in Science and Engineering XXI, vol.~98,
  pp. 467--475. Springer-Verlag (2014)

\bibitem{Mandal2}
Mandal, B.C.: Neumann-{N}eumann {W}aveform {R}elaxation {A}lgorithm in
  {M}ultiple {S}ubdomains for {H}yperbolic {P}roblems in 1d and 2d.
\newblock Numer. Methods Partial Differ. Equ.  (2016).
\newblock \urlprefix\url{DOI 10.1002/num.22112}

\bibitem{MR998911}
Marini, L.D., Quarteroni, A.: A relaxation procedure for domain decomposition
  methods using finite elements.
\newblock Numer. Math. \textbf{55}(5), 575--598 (1989).
\newblock \doi{10.1007/BF01398917}.
\newblock \urlprefix\url{http://dx.doi.org/10.1007/BF01398917}

\bibitem{MarQuar02}
Martini, L., Quarteroni, A.: An {I}terative {P}rocedure for {D}omain
  {D}ecomposition {M}ethods: a {F}inite {E}lement {A}pproach.
\newblock SIAM, in Domain Decomposition Methods for PDEs, I pp. 129--143 (1988)

\bibitem{pswr-dd22}
Ong, B.W., High, S., Kwok, F.: Pipeline schwarz waveform relaxation.
\newblock In: T.~Dickopf, M.J. Gander, L.~Halpern, R.~Krause, L.~Pavarino
  (eds.) Domain Decomposition Methods in Science and Engineering XXII, Lecture
  Notes in Computational Science and Engineering, pp. 179--187. Springer
  International Publishing (2016).
\newblock \doi{10.1007/978-3-319-18827-0_36}.
\newblock \urlprefix\url{http://dx.doi.org/10.1007/978-3-319-18827-0_36}

\bibitem{Picard}
Picard, E.: Sur l'application des m\'ethodes d'approximations successives \`a
  l'\'etude de certaines \'equations diff\'erentielles ordinaires.
\newblock Journal de Math\'ematiques Pures et Appliqu\'ees pp. 217--272 (1893)

\bibitem{MR2104179}
Toselli, A., Widlund, O.: Domain decomposition methods---algorithms and theory,
  \emph{Springer Series in Computational Mathematics}, vol.~34.
\newblock Springer-Verlag, Berlin (2005).
\newblock \doi{10.1007/b137868}.
\newblock \urlprefix\url{http://dx.doi.org/10.1007/b137868}

\bibitem{MR1340665}
Vandewalle, S.G., Van~de Velde, E.F.: Space-time concurrent multigrid waveform
  relaxation.
\newblock Ann. Numer. Math. \textbf{1}(1-4), 347--360 (1994).
\newblock Scientific computation and differential equations (Auckland, 1993)

\end{thebibliography}
\bibliographystyle{spmpsci}

\section*{Appendix}
\setcounter{table}{0}
\renewcommand{\thetable}{A\arabic{table}}
\begin{table}[htbp]
  \centering
  \begin{tabular}{|c|c|c|x{0.7in}|x{0.7in}|c|}
    \hline
    $J$ & Walltime (s) & Speedup & Actual Parallel Efficiency &
    Theoretical Peak Efficiency & MFLOPS/$\mu$s \\
    \hline
    8    & 492 & 6.82  & 0.43 & 0.53 & 46 \\
    16   & 402 & 8.35  & 0.52 & 0.70 & 57 \\
    32   & 357 & 9.40  & 0.59 & 0.82 & 64 \\
    64   & 233 & 14.40 & 0.90 & 0.90 & 98 \\
    128  & 221 & 15.18 & 0.95 & 0.95 & 103 \\
    256  & 216 & 15.53 & 0.97 & 0.98 & 105 \\
    512  & 213 & 15.75 & 0.98 & 0.99 & 107 \\
    1024 & 212 & 15.83 & 0.99 & 0.99 & 107 \\
    2048 & 211 & 15.90 & 0.99 & 1.00 & 108 \\
    4096 & 210 & 15.98 & 1.00 & 1.00 & 108 \\
    8192 & 211 & 15.90 & 0.99 & 1.00 & 108 \\
    \hline
  \end{tabular}
  \caption{Raw data used to compute the efficiency of the pipeline
    NNWR implementation as a function of the number of time blocks,
    $J$.  To compute the speedup and efficiency, the classical NNWR
    implementation with 8 processing cores was used -- this benchmark
    computation took 3355 seconds.}
  \label{tbl:efficiency_nn}
\end{table}

\begin{table}[htbp]
  \centering
  \begin{tabular}{|c|c|c|c|c|c|c|}
    \hline
    $K$ & \# procs & Walltime (s) & Parallel Efficiency &  MFLOPS/$\mu$s \\
    \hline
    1   & 16   & 211 & 1.00 & 27 \\
    2   & 32   & 209 & 1.01 & 54 \\
    4   & 64   & 214 & 0.99 & 108 \\
    8   & 128  & 212 & 1.00 & 213 \\
    16  & 256  & 212 & 1.00 & 429 \\
    32  & 512  & 212 & 1.00 & 815 \\
    64  & 1024 & 212 & 1.00 & 1540 \\
    \hline
  \end{tabular}
  \caption{Raw data used to compute the weak scaling capability of the
    pipeline NNWR implementation.  To compute the efficiency, the
    classical NNWR implementation with 8 processing cores compute one
    full iteration using 423 seconds.  Here, an efficiency of 1 means
    that the pipeline NNWR implementation is able to compute $K$ full
    waveform iterates using $2NK$ processing cores in half the time it
    takes $N$ processing cores to compute one full waveform iterate.}
    \label{tbl:weak_efficiency_nn}
\end{table}

\begin{algorithm}[htbp]
  \tcc{Assumes $J<2K$}
  \ForPar{$i=1$ \KwTo $N$}{
    $niter = \lceil 2K/J \rceil$\;
    \For{$iter=1$ \KwTo $niter$} {
      \ForPar{$l=1$ \KwTo $J$} {
        $a = (iter - 1)J + l$\;
        \If {$a \le 2K$} {
          $k = \lceil a/2 \rceil$\;
          $q = mod(a,2)$\;
          
          \If{$q==1$}{
            \tcc{Dirichlet Step: lines 5--19 in Algorithm~\ref{alg:pipeline_NNWR}}
          } \Else {
            \tcc{Auxiliary Step: lines 21--29 in Algorithm~\ref{alg:pipeline_NNWR}}
          }
        }
      }
    }
  }
  \caption{This pipeline implementation of the NNWR algorithm is able
    to utilize $NJ$ computing cores if the domain is broken into $N$
    non-overlapping subdomains and $K$ full iterates are used,
    provided $J<2K$.}
  \label{alg:pipeline_NNWR_case2}  
\end{algorithm}

\begin{algorithm}[htbp]
  \tcc{assumes $2K\le \lceil N/2 \rceil$}
  \KwIn{$N$: \# subdomains; $K$: \# waveform iterates}
  $m = \lceil N/2 \rceil$\;
  
  \ForPar{$p=1$ \KwTo $2K$}{
    
    $niter = \lceil N/(4K) \rceil$\;

    \For{$iter=1$ \KwTo $niter$} {

      \For{$k=1$ \KwTo $K$}{
        \If {$p > K$ } {
          $i = m + (iter-1)*2K + 2(p-K-1)$
        } \Else {
          $i = m + (iter-1)*2K + 2(p-K -1) + 1$
        }
        \If {$i\in[1,N]$} {
          \tcc{receive boundary info: lines 9--23 in Alg.~\ref{alg:DNWR_naive}}
          \tcc{Solve space time pde: lines 24--31 in Alg.~\ref{alg:DNWR_naive}}
          \tcc{Send boundary info: lines ~32--44 in Alg.~\ref{alg:DNWR_naive}}
        }

        \If {$p > K$ } {
          $i = m + (iter-1)*2K + 2(p-K-1) + 1$
        } \Else {
          $i = m + (iter-1)*2K + 2(p-K -1) $
        }
        \If {$i\in[1,N]$} {
          \tcc{receive boundary info: lines 9--23 in Alg.~\ref{alg:DNWR_naive}}
          \tcc{Solve space time pde: lines 24--31 in Alg.~\ref{alg:DNWR_naive}}
          \tcc{Send boundary info: lines ~32--44 in Alg.~\ref{alg:DNWR_naive}}
        }
        Check for convergence.  If converged, {\tt break}\; 
      }
    }
  }
  
  \caption{A classical implementation of the DNWR algorithm using $2K$
    processes, where $K$ is the number of waveform iterates being
    computed.  This pseudo code is for the case $2K\le \lceil N/2
    \rceil$.}
  \label{alg:DNWR_case2}  
\end{algorithm}

\begin{table}[htbp]
  \centering
  \begin{tabular}{|c|c|c|c|c|c|c|}
    \hline
    $K$ & \# procs & Walltime (s) &   MFLOPS/$\mu$s \\
    \hline
    1   & 8   & 211 &  3.4 \\
    2   & 16   & 211 &  5.0 \\
    4   & 32   & 211 &  8.4 \\
    8   & 64  & 211 &  15.1 \\
    16  & 128  & 211 &  28.5 \\
    32  & 256  & 212 &  54.8 \\
    64  & 512 & 215 &  106.5 \\
    128  & 1024 & 220  & 206.8 \\
    \hline
  \end{tabular}
  \caption{Raw data used to compute the weak scaling capability of the
    pipeline DNWR implementation.  The pipeline DNWR allows for
    multiple waveform iterates to be computed in approximately the
    same wall-clock time as eight processors computing a single
    waveform iterate.}
    \label{tbl:weak_efficiency_dnwr}
\end{table}

\end{document}